\documentclass[11pt, letterpaper, reqno]{amsart}

\usepackage[a4paper,
            top=30mm,bottom=30mm,
            left=35mm,right=35mm]{geometry}

\usepackage{amsmath, amsthm, amsfonts, amssymb, amscd, mathtools, bm} 
\usepackage{xcolor} 
\usepackage[unicode=true]{hyperref} 
\usepackage[shortlabels]{enumitem} 

\theoremstyle{plain}
\newtheorem{theorem}{Theorem}[section]

\newtheorem{proposition}[theorem]{Proposition}

\theoremstyle{definition}
\newtheorem{definition}[theorem]{Definition}

\theoremstyle{remark}

\numberwithin{equation}{section}

\newcommand{\rd}{\mathrm{d}}

\newcommand{\N}{\mathbb{N}}
\newcommand{\R}{\mathbb{R}}

\title[]{Strong convergence of sequences with \\ vanishing relative entropy}

\author[]{Nuno J. Alves}
\address[Nuno J. Alves]{
      University of Vienna, Faculty of Mathematics, Oskar-Morgenstern-Platz 1, 1090 Vienna, Austria.}
\email{nuno.januario.alves@univie.ac.at}

\author[]{Jakub Skrzeczkowski}
\address[Jakub Skrzeczkowski]{
      University of Oxford, Mathematical Institute, Woodstock Road, Oxford, OX2 6GG, United Kingdom.}
\email{jakub.skrzeczkowski@maths.ox.ac.uk}

\author[]{Athanasios E. Tzavaras}
\address[Athanasios E. Tzavaras]{
      King Abdullah University of Science and Technology, CEMSE Division, Thuwal 23955-6900, Saudi Arabia.}
\email{athanasios.tzavaras@kaust.edu.sa}

\begin{document}

\begin{abstract}
We show that under natural growth conditions on the entropy function, convergence in relative entropy implies $L^p$-convergence; for $p>1$ this implication becomes an equivalence. The main tool is the theory of Young measures, in a form that accounts for the formation of concentrations in weak limits. 
\end{abstract}

\maketitle
\thispagestyle{empty} 

\section{Introduction}
For a continuously differentiable function $f : \mathbb{R}^n \to \mathbb{R}$, define $f( \cdot | \cdot )$ by the quadratic part of the Taylor series expansion
\begin{equation} \label{defrelfun}
f(v | u) = f(v) - f(u) - Df(u) \cdot (v - u)
\end{equation}
for $v, u \in \mathbb{R}^n$, and note that if $f$ is strictly convex, then $f(v|u) > 0$ for $v \ne u$. When $f$ is the entropy (energy),  the function $f( \cdot |  \cdot)$ is called \textit{relative entropy} (\textit{energy}). 
In several problems in mechanics,  kinetic theory, or information theory, the function $f$ is convex and one is able to derive information on the convergence of a sequence $\{ u_n \}$
to a limit function $u$ via an integral relation of the form
\begin{equation}\label{convrelen}
F [ u_n | u ] = \int_\Omega f ( u_n  | u  ) \, \rd x \to 0 \qquad \text{as} \ n \to \infty
\end{equation}
where $\Omega \subseteq \R^d$, $d \in \N$, is an open set (see {\it e.g.} \cite{dafermos79,gilardoni2010pinsker,unterreiter2000generalized}). The question arises as
to what type of convergence this information translates into. An attempt in this direction, employing real analysis techniques, may be found in \cite{alves2024mode}.

The reader should note that $f(\cdot | \cdot )$ or its integral does not define a distance and is not even symmetric. The lack of symmetry may be alleviated by introducing the
symmetrized form of the relative entropy:
$$
f_{\text{sym}} (v|u) = f(v | u) + f(u | v) = ( Df(u) - Df(v) ) \cdot (u - v)
$$
which is symmetric and (strictly) positive for $f$ strictly convex and $v \ne u$. However, $(v,u) \mapsto \int f_{\text{sym}}(v  |u  ) \, \rd x$ still does not induce a distance, as it does not satisfy the triangle inequality.
The intent of this note is to prove that under natural growth conditions of the entropy $f$, convergence in \eqref{convrelen} implies strong convergence in $L^p$.

This result can be easily checked in some special cases. It is trivial in the quadratic case $f(u) = |u|^2$, since then
 $f(v | u) = |v - u|^2$. Another interesting case concerns $h(u) = u \ln u$ (for $n = 1$). The relative entropy then takes the form
 \[
 h(v | u) = v \ln \frac{v}{u} - (v - u) \, .
 \]
 When $v, u$ are probability densities, $\int v \, \rd x = \int u \, \rd x = 1$, the relevant quantity becomes the functional 
 \begin{equation} \label{rel_entropy_H}
 H ( v | u ) = \int v \ln \frac{v}{u} \, \rd x
 \end{equation}
 which is often used in information theory and statistics under the name Kullback--Leibler divergence. The Csisz{\'a}r--Kullback--Pinsker inequality states that for $v,u \in L^1(\R^d)$ nonnegative functions with
 $\| v \|_1 = \| u \|_1 = 1$ it holds that (see \cite{gilardoni2010pinsker})
 \begin{equation}\label{CKP_inequality}
 \| v - u \|_1 ^2 \le 2 H(v | u)
 \end{equation}
 indicating that convergence in relative entropy implies $L^1$-convergence. \par Inequality \eqref{CKP_inequality} can be generalized as 
 $$
 \frac{A}{2^{2/p}} \, \min(\|v\|_{p}^{p-2}, \|u\|_{p}^{p-2}) \, \|u-v\|^2_{p} \leq \int f(v|u) \, \rd x
 $$
 where $u, v \geq 0$, $u, v \in L^1 \cap L^p $, $p \in [1,2]$ and $f \in C^2(0, \infty)$ is a strictly convex nonnegative function with $f(1)=f'(1) = 0$ and $ A := \inf_{s\in(0,\infty)} s^{2-p}\,f''(s)>0$; see \cite[Proposition 3.1]{caceres2002nonlinear}.
 \par
Our approach employs the theory of Young measures and a representation theorem from~\cite{demoulini2012weak}
as methodological tools to address this question.
The analysis shows that convergence in relative entropy precludes concentrations in weak limits. \par 
The manuscript is organized as follows. Section~\ref{section_main_results} contains the statement of the main results. The first result establishes a representation formula for convex entropies that captures potential concentrations in the limit of bounded sequences; its proof is given in Section~\ref{section_proof_1}. The second and third results relate convergence in relative entropy to strong convergence in $L^p$, and their proofs are presented in Section~\ref{section_proof_2} and Section~\ref{section_proof_3}, respectively.

\section{Preliminaries and statement of the main results} \label{section_main_results}

We start by presenting two situations where one may conclude immediately. Suppose the function $f : \mathbb{R}^n \to \mathbb{R}$ satisfies the
hypothesis that its Hessian is bounded from below
\begin{equation}\label{unifconv}
D^2 f (\lambda) \ge c I  
\end{equation}
for some $c > 0$ and all $\lambda \in  \mathbb{R}^n$, and let $\{ u_n \}$ be a sequence satisfying \eqref{convrelen}. We then conclude
\[
c \int_\Omega |u_n - u|^2 \, \rd x \le  \int_\Omega f( u_n | u ) \, \rd x  \to 0 \quad \text{as} \  n \to \infty.
\]
\par 
Another direct case occurs when $f(\lambda) = |\lambda|^p$, $p > 1$,  and we assume for the sequence $\{ u_n \}$ that
\begin{equation}\label{hypo2}
u_n \rightharpoonup u \quad \text{weakly in} \ L^p(\Omega) \quad \text{and} \quad f(u_n | u) \to 0 \quad \text{in} \ L^1(\Omega)
\end{equation}
for some $u \in L^p(\Omega)$. Then, from \eqref{defrelfun} and the hypothesis \eqref{hypo2} we deduce
\[
\int_\Omega |u_n|^p \, \rd x - \int_\Omega |u|^p \, \rd x = \int_\Omega p |u|^{p-2} u \cdot (u_n - u) \, \rd x + \int_\Omega f(u_n | u) \, \rd x  \to 0
\]
that is, $ \| u_n \|_{L^p(\Omega)} \to \| u \|_{L^p(\Omega)}$. Hence, $u_n \to u$ in $L^p(\Omega)$.

The assumption of global uniform convexity \eqref{unifconv} might be too restrictive (and often unattainable when $L^\infty$ bounds
are not available), while the case of the $L^p$ norm is specialized. To address the problem for more general functions, we use the theory of Young measures introduced by L. Tartar, who formalized an idea for generalized functions initially proposed by L. C. Young. For thorough presentations of the relevant results, the reader is referred to {\sc Ball} \cite{ball1989version} and {\sc 
Pedregal} \cite[Chapter 6]{pedregal1997parametrized}. The second key 
component is the representation of potential concentrations in the 
entropy norm via a concentration measure, as developed in~\cite[Appendix A]{demoulini2012weak}. This representation is generalized here to hold for a convex entropy $h$.

\begin{theorem}\label{prop1}
Let $h: \mathbb{R}^n \to \R$ be a convex function that
satisfies for $p \geq 1$ the bound
\begin{equation*} \label{growth0}
h(\lambda) \geq \frac{1}{c}|\lambda|^p - c 
\end{equation*}
for some constant $c > 0$ and every $\lambda \in \mathbb{R}^n$. Let $\{ u_n \}$ be a sequence of functions $u_n : \Omega \to \R^n$,
where $\Omega$ is an open bounded subset of $\R^d$, that satisfy the uniform bounds
\begin{equation}\label{unibound}
\sup_n \int_\Omega h(u_n) \, \rd x  \le  C  < \infty \, .
\end{equation}
There exist a subsequence $\{ u_{n_k} \}$, a parametrized family of probability measures (Young measures) $\{ \nu_x \}_{x \in \Omega} $,
and a concentration measure $\mu (\rd x)$ such that
\begin{equation}\label{ymrep}
 f(u_{n_k})  \rightharpoonup   \int f(\lambda) \, \rd \nu_x(\lambda)   \ \mbox{ weakly in $L^1(\Omega)$}
\end{equation}
for every continuous function $f : \R^n \to \R$ satisfying $\lim_{|\lambda| \to \infty} \frac{|f(\lambda)|}{1+ |\lambda|^p} = 0$, and
\begin{equation}\label{concmeas}
 \int_\Omega h(u_{n_k}) \, \varphi \, \rd x  \to    \int_\Omega \Big ( \int h(\lambda) \, \rd \nu_x(\lambda) \Big ) \, \varphi \, \rd x +  \int_\Omega \varphi \, \mu (\rd x) \quad \mbox{as $k \to \infty$}
\end{equation}
for all $\varphi \in C(\bar{\Omega})$.
\end{theorem}
 
We remark that the representation \eqref{ymrep} of the weak limit via Young measures in the most general setting holds for sequences $\{u_n\}$ and $\{\psi(x,u_n)\}$ such that $\{u_n \}$ is bounded in $L^q(\Omega)$ for some $q > 0$, and $\psi$ is a Carath\'{e}odory function so that $\{\psi(x, u_n(x))\}$ is relatively weakly compact in $L^1(\Omega)$; see \cite[Theorem 6.2]{pedregal1997parametrized}. In that case, we have
\[
\psi(x,u_{n_k}(x))  \rightharpoonup   \int \psi(x,\lambda) \, \rd \nu_x(\lambda)   \ \mbox{ weakly in } L^1(\Omega).
\]
\par 
Theorem~\ref{prop1}, combined with the characterizations of uniform $p$-integrability for a sequence of functions $\{u_n \}$ (see Appendix~\ref{appendix_B}), leads to an equivalence between $L^p$-convergence for $p>1$ and convergence in relative entropy. 
 
\begin{theorem} \label{thmp>1}
Let $h: \mathbb{R}^n \to \R$ be a continuously differentiable, strictly convex function satisfying the bounds
\begin{equation} \label{growth}
 \frac{1}{c}|\lambda|^p - c \leq h(\lambda) \leq c |\lambda|^p + c
\end{equation}
for some $p > 1$, $c>0$ and every $\lambda \in \R^n$. \par 
Let $u_n, u : \Omega \to \mathbb{R}^n$ belong to $L^p(\Omega)$, where $\Omega$ is an open bounded subset of $\R^d$. Then 
\begin{equation*}
\int_\Omega h(u_n | u) \, \rd x \to 0 \qquad \text{if and only if} \qquad u_n \to u \ \text{in} \ L^p(\Omega)
\end{equation*}
as $n \to \infty$.
\end{theorem}
To the best of our knowledge, Theorem~\ref{thmp>1} has not appeared in the literature in this precise form, although it may be regarded as folklore. After establishing it as an application of Theorem~\ref{prop1}, we later learned of a different proof, kindly communicated to us by a referee. For completeness, and since it highlights a distinct line of reasoning, we include this alternative argument in Appendix~\ref{appendix_A}.
 \par 
An obvious example satisfying the conditions of Theorem \ref{thmp>1} is the function $h(\lambda) =  |\lambda|^p$. Convergence in the corresponding relative entropy is precisely the type of control that arises in several relaxation and singular limit problems for Euler-type systems; see for instance~\cite{lattanzio2017gas, alves2024zero, alves2024weak}. 
\par
On the other hand, in applications of relative entropy, such as those arising in kinetic theory or information theory, one often encounters the function
$\lambda \ln \lambda$. The next result covers this situation by considering the shifted function
\[h(\lambda) = 1/e + \lambda \ln \lambda, \quad \lambda \geq 0,\]
and is derived using Theorem~\ref{prop1}.

\begin{theorem} \label{thmp=1}
Consider a function $h : [0,\infty) \to [0, \infty)$ and assume the following: 
\begin{enumerate}[(i)]
\item $h \in C([0,\infty)) \cap C^2(0,\infty)$,
\vspace{2mm}
\item $h^{\prime \prime}(\lambda) > 0$ for every $\lambda > 0$,
\vspace{2mm}
\item $ \dfrac{h(\lambda)}{\lambda} \to \infty$ as $\lambda \to \infty$,
\vspace{2mm} 
\item $c\lambda - c \leq h(\lambda)$ for some $c > 0$ and every $\lambda \geq 0$.
\end{enumerate}
Let $\{ u_n \}\subseteq L^1(\Omega)$ and $u \in L^\infty(\Omega)$ be nonnegative. If $Dh(u) \in L^\infty(\Omega)$,
\begin{equation} \label{hbound}
\sup_n \int_\Omega h(u_n) \, \rd x = C < \infty
\end{equation}
and 
\[\int_\Omega h(u_n | u) \, \rd x \to 0 \qquad \text{as} \ n \to \infty\]
then $\{ u_n\}$ converges to $u$ in $L^1(\Omega)$.
\end{theorem}

\section{Representation via Young and concentration measures} \label{section_proof_1}

We provide the proof of Theorem~\ref{prop1}.  Let $\{ u_n \}$ satisfy the uniform bounds \eqref{unibound}. The standard theory of Young measures 
\cite{ball1989version}, \cite{pedregal1997parametrized}  implies there exists a
subsequence $\{ u_{n_k} \}$ and  a parametrized family of probability measures $\{ \nu_x \}_{x \in \Omega}$ such that
\eqref{ymrep} holds for any continuous $f : \R^n \to \R$ that satisfies $\lim_{|\lambda| \to \infty} \frac{|f(\lambda)|}{1+ |\lambda|^p} = 0$. \par
The representation formula does not cover the entropy function $h$. There are two alternatives that one may pursue. One can use the theory
of concentration measures of {\sc DiPerna-Majda} \cite{dipm87} or {\sc Alibert-Bouchitt{\'e}} \cite{ab97}. This approach has the advantage that it can represent entire families of functions  with critical growth but under the assumption that their recession functions satisfy continuity properties on the unit sphere. 
An alternative approach is pursued in \cite{demoulini2012weak} which can only represent the entropy function, but with no continuity 
assumptions on the recession function, only assuming convexity and positivity. 
We assume first that $h : \R^n \to \R$ is convex and nonnegative (the nonnegativity assumption will be removed later).
The proof  proceeds in two steps:

First, one considers the map $\langle \nu , h \rangle : x \mapsto \int h(\lambda) \, \rd \nu_x (\lambda)$ for $x \in \Omega$ and shows that 
\begin{equation}\label{belongsL1}
\int h(\lambda) \, \rd \nu_x (\lambda) \in L^1 (\Omega) \, .
\end{equation}
To this end, for $R>0$, let $h_R$ be the truncation of $h$:
\[
h_R (\lambda) = 
\begin{cases} 
h(\lambda)  &\mbox{when} \quad h(\lambda) \leq R,
\\
R           &\mbox{when} \quad h(\lambda) >  R.
 \end{cases}
 \]
Since $h_R  \nearrow h$ as $R \to \infty$, the monotone convergence theorem implies 
$$
\int h(\lambda) \, \rd \nu_x (\lambda) := \lim_{R \to \infty} \int h_R (\lambda) \, \rd \nu_x (\lambda)  \quad \mbox{for a.e. $x \in \Omega$}
$$
providing a definition for the bracket $\langle \nu , h \rangle$. Observe that
$$
\begin{aligned}
\int_\Omega \Big ( \int h_R (\lambda) \, \rd \nu_x (\lambda) \Big ) \, \rd x &= \lim_{k \to \infty} \int_\Omega h_R (u_{n_k}) \, \rd x
\\
&\le \limsup_{k \to \infty} \int_\Omega h(u_{n_k}) \, \rd x \le C \, .
\end{aligned}
$$
Letting $R \to \infty$ and using $h_R  \nearrow h$  yields \eqref{belongsL1}.

Second, for the sequence $\big \{ h(u_{n_k}) - \int h(\lambda) \, \rd \nu_x (\lambda) \big \}$ we prove consecutively:
\begin{itemize}
\item[(i)]   $\sup_k   \int_\Omega \left | h(u_{n_k}) - \int h(\lambda) \, \rd \nu_x (\lambda) \right | \, \rd x  \le C < \infty $,
\item[(ii)] There exists a concentration measure $\mu(\rd x) \in \mathcal{M}(\Omega)$ such that along a subsequence
\[
h(u_{n_k}) - \int h(\lambda) \, \rd \nu_x (\lambda) \rightharpoonup \mu(\rd x) \qquad \mbox{ weak-$\ast$ in measures,}
\]

\item[(iii)] $\mu (\rd x ) \ge 0$.
\end{itemize}

Part (i) follows from \eqref{belongsL1} and \eqref{unibound}. Then, due to (i), there exists a (signed) measure 
$\mu \in \mathcal{M}(\Omega)$ such that
along a subsequence which we still denote by $\{ u_{n_k} \}$ we have
$$
\int_\Omega \left ( h(u_{n_k}) - \int h(\lambda) \, \rd \nu_x (\lambda) \right ) \varphi(x) \, \rd x  \to \int_\Omega \varphi(x) \, \mu (\rd x )
$$
for any $\varphi \in C(\bar{\Omega})$.
Finally, for $\varphi \in C(\bar{\Omega})$, $\varphi \ge 0$, we have
$$
\begin{aligned}
\int_\Omega \left (  \int h(\lambda) \, \rd \nu_x (\lambda) \right ) \varphi(x) \, \rd x  
&= \sup_{R > 0} \int_\Omega \left (  \int h_R (\lambda) \, \rd  \nu_x (\lambda) \right ) \varphi(x) \, \rd x  
\\
&= \sup_{R > 0} \lim_{k \to \infty} \int_\Omega h_R(u_{n_k}) \, \varphi (x) \, \rd x 
\\
&\le \lim_{k \to \infty} \int_\Omega h(u_{n_k}) \,  \varphi (x) \, \rd x 
\end{aligned}
$$
which implies $\mu := \mbox{weak}-\ast \lim \left ( h(u_{n_k}) - \int h(\lambda) d \nu_x (\lambda) \right ) \ge 0$.

This completes the proof when $h$ is a convex and nonnegative function.  
If we only know that $h$ is convex, then  there exists $a \in \R^n$, $b \in \R$ such that
$$
h(\lambda) \ge - a\cdot \lambda - b \, .
$$
The function $\hat h(\lambda) := h(\lambda) + a \cdot \lambda + b$ is convex and nonnegative. We apply the previous procedure to  
$\hat h(\lambda)$. Since $u_{n_k} \rightharpoonup u$ weakly in $L^1(\Omega)$ we can transfer the relation \eqref{concmeas} from the function
$\hat h$ to the function $h$ and complete the proof.

\section{Proof of Theorem~\ref{thmp>1}} \label{section_proof_2}
First, we prove that strong convergence in $L^p(\Omega)$ implies convergence in relative entropy. Suppose that $\{ u_n \}$ converges to $u$ in $L^p(\Omega)$. From  \cite[Proposition 2.32]{dacorogna2008direct}  and (\ref{growth}), it follows that $Dh(u) \in L^{p^\prime}(\Omega)$, where $p^\prime = \frac{p}{p-1}$ the  conjugate exponent of $p$. Hence
\begin{equation} \label{aux1}
\int_\Omega |Dh (u)| |u_n - u| \, \rd x  \to 0 \quad \text{as} \ n \to \infty.
\end{equation}
Using \cite[Proposition 2.32]{dacorogna2008direct} once more, 
\begin{equation} \label{aux2}
\begin{split}
\int_\Omega  |h (u_n)  - h (u)| \, \rd x 
 &\leq  C \int_\Omega (1 + |u_n|^{p-1} + |u|^{p-1}) |u_n - u| \, \rd x \\
& \leq C \left( \int_\Omega (1 + |u_n|^{p-1} + |u|^{p-1})^{p^\prime} \, \rd x \right)^{\frac{1}{p^\prime}} \|u_n - u \|_{L^p(\Omega)} \\
& \to 0 \quad \text{as} \ n \to \infty.
\end{split}
\end{equation}
Combining (\ref{aux1}) and (\ref{aux2}) yields the desired convergence in relative entropy.

The converse direction is more intricate. We prove that every subsequence of $\{u_n \}$ has a further subsequence that converges in $L^p(\Omega)$ to $u$. Thus, it suffices to show the existence of a subsequence of $\{u_n\}$ that converges to $u$ in $L^p(\Omega).$ Let $\{ u_n \}$, $u$ satisfy
\begin{equation}\label{relenconv}
\int_\Omega h(u_n | u ) \, \rd x \to 0 \, .
\end{equation}
Assume without loss of generality that $h(\lambda) \ge 0$. Since $h$ is convex, it satisfies for some
$a \in \R^n$, $b \in \R$ the lower bound
$$
h (\lambda) \ge - a\cdot \lambda - b.
$$
If $h$ assumes negative values, then it is replaced by the convex, nonnegative function 
$\hat h(\lambda) = h (\lambda) + a\cdot \lambda + b $ which still satisfies \eqref{relenconv}. \par 
We first prove that $\int_{\Omega} h(u_n) \, \rd x \leq C$ for some $C > 0$ and every $n \in \N$ .  From the convergence in relative entropy, the sequence $\left\{\int_\Omega h(u_n | u) \, \rd x \right\}$ is bounded, so there exists $K > 0$ such that for every $n \in \N$,
\begin{equation}\label{aux3}
\begin{split}
     \int_\Omega h(u_n) \, \rd x  = & \,\int_\Omega h(u_n | u) + h(u) + Dh(u) \cdot (u_n - u) \, \rd x \\
     \leq & \  K + c \int_\Omega |u|^p \, \rd x + c|\Omega| + \int_\Omega |Dh(u) \cdot u| \, \rd x + \int_\Omega |Dh(u) \cdot u_n| \, \rd x \, .
\end{split}
\end{equation} 
We estimate the last term on the right-hand side of the expression above using Young's inequality with a parameter $\epsilon$ to obtain
\begin{equation} \label{aux4}
\int_\Omega |Dh(u) \cdot u_n| \, \mathrm{d}x \leq C(\epsilon)\int_\Omega |Dh(u)|^{p^\prime} \mathrm{d}x + \epsilon c^2 |\Omega| + \epsilon c \int_\Omega h(u_n) \, \mathrm{d}x 
\end{equation}
where $C(\epsilon) = (\epsilon p)^{-\frac{p^\prime}{p}} (p^\prime)^{-1}$. Combining (\ref{aux3}) with (\ref{aux4}), taking $\epsilon = 1/(2c)$ and noting that $Dh(u)\cdot u \in L^1(\Omega)$ yields the desired uniform bound. \par
Let $\{ u_{n_k} \}$ be a subsequence, $\{ \nu_x \}$ parametrized probability measures and $\mu(\rd x)$ a concentration measure
as in Theorem~\ref{prop1}.  Observe that $u_{n_k} \rightharpoonup \int \lambda \, \rd \nu_x (\lambda)$ weakly in $L^p(\Omega)$, and hence
$$
\lim_{k \to \infty}\int_{\Omega} Dh(u)\cdot (u_{n_k} - u) \, \rd x = \int_{\Omega} \int Dh(u)\cdot (\lambda - u) \, \rd \nu_x (\lambda) \,\rd x.
$$
Therefore, using \eqref{concmeas} with $\varphi = 1$, we have
$$
\int_{\Omega} \int h \big (\lambda | u \big ) \, \rd \nu_x (\lambda) \, \rd x  +  \mu(\Omega) = 0
$$
which implies $\mu = 0$ and $\int h(\lambda | u) \, \rd \nu_x (\lambda) = 0$ for a.e. $x \in \Omega$ since $h (\lambda | u) \geq 0$. 
The strict convexity of $h$ then yields $\nu_x  = \delta_{u(x)}$. Applying \eqref{ymrep} to the function $\psi(\lambda, x) = | \lambda -  u(x) |^q$,
with $1\le q < p$, gives 
$$
\int_\Omega  | u_{n_k}  - u  |^q \, \rd x \to 0 
$$
and along subsequences $u_{n_k} \to u $ almost everywhere and in measure. Since all subsequences of $\{ u_n \}$ have the same limit we
conclude $u_n \to u$ in $L^q(\Omega)$, $q \in [1, p)$. \par

The final step is to establish the convergence in $L^p (\Omega)$ of a subsequence of $\{u_{n_k} \}$. For that subsequence (not relabeled), it suffices to prove that $\{u_{n_k} \}$ is uniformly $p$-integrable; see the Appendix for characterizations of this property. The result then follows by the Vitali convergence theorem. \par  
For simplicity of notation, we relabel the sequence $\{u_{n_k} \}$ as $\{u_k \}$ for the remainder of the proof. Since  $u_k \rightharpoonup u$ weakly in $L^p(\Omega)$ it follows that
\begin{equation} 
\int_\Omega Dh (u) \cdot (u_k - u) \, \rd x  \to 0 \quad \text{as} \ n \to \infty.
\end{equation}
In turn, \eqref{relenconv} implies
\begin{equation}\label{aux5}
\int_\Omega h(u_{k} ) \, \rd x  \to \int_\Omega h(u)\, \rd x  \, .
\end{equation}
Now, since $\{u_k \}$ converges to $u$ in measure, then along a subsequence if necessary, $\{u_k \}$ converges to $u$ almost everywhere. For each $M > 0$, define $\varphi_M : \mathbb{R}^n \to \mathbb{R}$ by
\begin{equation}
\varphi_M(\lambda) = 
\begin{cases}
1, \quad & \text{if} \ 0 \leq |\lambda| \leq M, \\
M - |\lambda|+1, \quad & \text{if} \ M < |\lambda| \leq M + 1, \\
0, \quad & \text{if} \  |\lambda| \geq M + 1. 
\end{cases}
\end{equation}
The sequence of functions  $\{ h(u_k) \, \varphi_M(u_k) \}$ is dominated by the constant function 
$\max_{0 \leq |\lambda| \leq M+1} h(\lambda)$. By the dominated convergence theorem, 
\begin{equation*} 
\int_\Omega h(u_k) \varphi_M(u_k) \, \rd x  \to \int_\Omega h(u) \varphi_M(u) \, \rd x  \quad \text{as} \ k \to \infty, 
\end{equation*}
which together with (\ref{aux5}) implies
\begin{equation*}
\int_\Omega h(u_k)(1- \varphi_M(u_k)) \, \rd x  \to \int_\Omega h(u) (1-\varphi_M(u)) \, \rd x  \quad \text{as} \ k \to \infty.
\end{equation*}
We conclude that for each $M > 0$ there exists $K_M \in \mathbb{N}$ such that for $k \geq K_M$ it holds 
 \begin{equation} \label{aux6}
\int_\Omega h(u_k) (1 - \varphi_M(u_k)) \, \rd x  \leq 2 \int_\Omega h(u)(1 - \varphi_M(u)) \, \rd x . 
 \end{equation}
Next, using $(\ref{growth})$, we select $M_\ast > 1$ such that 
$$
h (\lambda)  \geq \frac{1}{2c} |\lambda|^p \quad \mbox{for $|\lambda| \geq M_\ast$ } . 
$$
For $M > M_\ast$ we may select $K_M$ such that for $k > K_M$ using \eqref{aux6} we have 
\begin{align} 
\int_{ \{ |u_k| \geq M + 1 \} } |u_k|^p \, \rd x  & \leq 2 c  \int_{ \{ |u_k| \geq M + 1 \} } h(u_k) \, \rd x  
\nonumber
\\
& =  2c  \int_{ \{ |u_k| \geq M + 1 \} } h(u_k)(1 - \varphi_M(u_k)) \, \rd x  
\nonumber
\\
& \leq 4 c  \int_\Omega h(u)(1 - \varphi_M(u)) \, \rd x . \label{aux10}
\end{align}
Now, let $\varepsilon > 0$ be fixed and select $\tilde M = \tilde{M} (\varepsilon)  > M_\ast$ so that 
\begin{equation*} 
4 c  \int_\Omega  h(u)(1 - \varphi_{\tilde{M}}(u)) \, \rd x  < \varepsilon^p. 
\end{equation*}
For $k \geq K_{\tilde{M}}$, \eqref{aux10} gives 
\begin{equation}\label{aux7}
\int_{ \{ |u_k| \geq \tilde{M}+ 1 \} } |u_k|^p \, \rd x   <\varepsilon^p.
\end{equation}
Since $u_k \in L^p(\Omega)$, 
for each  $ k =1, \ldots, K_{\tilde{M}} - 1$ there exists $M_k > 0$ so that
\begin{equation} \label{aux8}
\int_{ \{ |u_k| \geq M_k \} } |u_k|^p \, \rd x   \leq \varepsilon^p  \, ,  \quad k = 1, \ldots , K_{\tilde{M}} - 1 \, .
\end{equation}
Setting $M = \max\{ \tilde{M}+ 1, M_1 , \ldots, M_{K_{\tilde{M} }} - 1 \}$, from (\ref{aux7}) and (\ref{aux8}) it follows that 
\begin{equation*}
\sup_k \int_{ \{ |u_k| \geq M \} } |u_k|^p \, \rd x   < \varepsilon^p 
\end{equation*}
which proves that $\{ u_k \}$ is uniformly $p$-integrable. This completes the proof of Theorem \ref{thmp>1}.

\section{Proof of Theorem~\ref{thmp=1}} \label{section_proof_3}
Similarly to the proof of Theorem \ref{thmp>1}, we show that $\{ u_n \}$ has a subsequence that converges to $u$ in $L^1(\Omega)$. For simplicity, none of the subsequences will be relabeled. We remark that we are working here with a nonnegative sequence $u_n \geq 0$.  \par 
First, we note that the sequence $\{ u_n \}$ is bounded in $L^1(\Omega)$. This follows from condition $(iv)$ in conjunction with (\ref{hbound}). Next, we prove that $\{ u_n \}$ is uniformly integrable. This follows from the De La Vall\'{e}-Poussin criterion;
 we provide the details for the reader's convenience. Let $\varepsilon > 0$ and choose $M_\varepsilon > 0$ such that $C/M_\varepsilon < \varepsilon$. From $(iii)$ we infer the existence of a constant $\Lambda_\varepsilon > 0$ such that whenever $\lambda > \Lambda_\varepsilon$ one has $h(\lambda) \geq M_\varepsilon \lambda$. Therefore 
\begin{align*}
\int_{ \{u_n \geq \Lambda_\varepsilon \}} u_n \, \rd x  & \leq \frac{1}{M_\varepsilon}\int_{ \{u_n \geq \Lambda_\varepsilon \}} h(u_n) \, \rd x  \\
& \leq \frac{C}{M_\varepsilon} < \varepsilon 
\end{align*}
which establishes the uniform integrability of $\{ u_n \}$. Thus, along a subsequence if necessary, 
$\{ u_n \}$ converges weakly in $L^1(\Omega)$. There exists a Young measure $\{ \nu_x \}_{x\in \Omega}$ generated by a subsequence of $\{ u_n \}$, such that for any  Carath\'{e}odory function $\psi : \Omega \times \mathbb{R}^+ \to \mathbb{R} $, if  $\{ \psi(x,u_n(x)) \}$ converges weakly in $L^1(\Omega)$,
its weak limit is represented by the map
\[x \mapsto \int_{\mathbb{R}^+} \psi(x,\lambda) \, \rd  \nu_x(\lambda). \]
Noting that the arguments in the proof of Theorem~\ref{prop1} also work in this case with $p=1$, there exists a (nonnegative) concentration measure $\mu(\rd x)$ such that  
\begin{equation}\label{concmeas2}
 \int_\Omega h(u_{n}) \, \varphi \, \rd x  \to    \int_\Omega \Big ( \int_{\R^+} h(\lambda) \, \rd \nu_x(\lambda) \Big ) \, \varphi \, \rd x +  \int_\Omega \varphi \, \mu (\rd x)
\end{equation}
as $n \to \infty$, for all $\varphi \in C(\bar{\Omega})$. As in the proof of Theorem \ref{thmp>1} we deduce that $\nu_x= \delta_{u(x)}$ for almost every $x \in \Omega$. \par 
Now, we apply \cite[Proposition 6.5]{pedregal1997parametrized} to the Carath\'{e}odory function $\psi(x, \lambda) = |\lambda - u(x)|$. 
Since $u \in L^\infty(\Omega)$,  we have $\psi(x, \lambda) \leq K(1 + \lambda)$ for some positive constant $K$. 
We claim that the sequence $\{ |u_n - u| \}$ is relatively weakly compact in $L^1(\Omega)$. Clearly, it is bounded in $L^1(\Omega)$. 
Moreover, for each $k > K$, let $m_k = \frac{k}{K} - 1$ so that 
\[ \{ \lambda \in \mathbb{R}^+ \ : \ K(1 + \lambda) \geq k \} \subseteq \{ \lambda \in \mathbb{R}^+ \ : \ \lambda \geq m_k  \}. \]
From condition $(iii)$ it follows that 
\[ \lim_{\lambda \to \infty} \frac{h(\lambda)}{K(1 + \lambda)} = \infty \] 
and so we can find $M_k>0$, with $M_k  \to \infty$ as $k \to \infty$,  such that $h(\lambda) \geq M_k K(1 + \lambda)$ whenever $\lambda \geq m_k$. Consequently, 
\begin{align*}
\int_{ \{ \psi(x, u_n(x)) \geq k \} } \psi(x, u_n(x)) \, \rd x  & \leq \int_{ \{ \psi(x, u_n(x)) \geq k \} } K(1 + u_n) \, \rd x  \\
& \leq \int_{ \{  u_n  \geq m_k \} } K(1 + u_n) \, \rd x  \\
& \leq \frac{1}{M_k}\int_{ \{  u_n  \geq m_k \} } h(u_n) \, \rd x  \\
& \leq \frac{C}{M_k} \to 0 \qquad \text{as} \ k \to \infty 
\end{align*}
which establishes the claim. Hence, along a subsequence if necessary, $\{ |u_n - u| \}$ is weakly convergent in $L^1(\Omega)$, which in particular implies that
\begin{align*}
\int_\Omega |u_n - u| \, \rd x   & \to \int_\Omega \int_{\R^+} |\lambda - u| \, \rd \delta_{u(x)}(\lambda) \, \rd x  = 0  \qquad \text{as} \ n\to \infty
\end{align*}
and concludes the proof.

\section*{Acknowledgments}
We thank the anonymous referee for communicating to us the proof presented in Appendix~\ref{appendix_A}.
\appendix

\section{Alternative proof to Theorem~\ref{thmp>1}} \label{appendix_A}

Assuming that $\{u_n\}$ converges to $u$ in relative entropy, we show that $\{u_n\}$ converges to $u$ in measure and that $\{u_n \}$ is uniformly $p$-integrable; the conclusion then follows from the Vitali convergence theorem.
\par We first observe that $\{ u_n\}$ is bounded in $L^p(\Omega)$; this follows from the same reasoning as in \eqref{aux3} and \eqref{aux4}. We proceed to prove that $\{u_n \}$ is uniformly $p$-integrable. Let $E$ be any measurable subset of $\Omega$. We have 
\begin{equation}
\int_E |u_n|^p \, \rd x \leq c \int_E h(u_n | u) \, \rd x + c \int_E h(u) \, \rd x + cB \|\nabla h(u) \chi_E \|_{L^{p^\prime}(\Omega)} + c^2|E|
\end{equation}
where $B = \|u \|_{L^p(\Omega)} + \sup_n \|u_n \|_{L^p(\Omega)}$. Note that the sequence $\{h(u_n | u) \}$ is uniformly integrable (since it is convergent in $L^1(\Omega)$), the function $h(u)$ has absolutely continuous integrals (since it belongs to $L^1(\Omega)$), and $\nabla h(u)$ has absolutely continuous $p^\prime$-integrals. Let $\varepsilon > 0$. From what was just observed, there exists $\tilde \delta_\varepsilon > 0$ such that if $|E| < \tilde \delta_\varepsilon$, then
\[\sup_n \int_E h(u_n | u) \, \rd x < \frac{\varepsilon^p}{4c}, \qquad \int_E h(u) \rd x < \frac{\varepsilon^p}{4c} \]
and 
\[\| \nabla h(u) \chi_E \|_{L^{p^\prime}(\Omega)} < \frac{\varepsilon^p}{4cB}. \]
Consequently, for $|E| < \delta_\varepsilon = \min\{ \tilde \delta_\varepsilon, \varepsilon^p/(4c^2) \}$ we have 
\[\sup_n \int_E |u_n|^p \, \rd x < \varepsilon^p \]
which establishes the uniform $p$-integrability of $\{u_n \}$. \par 
Now we claim that for every $\varepsilon>0$ and every $K > 1$, there exists $\delta > 0$ such that for all $x,y \in \R^m$, if $|x-y| \geq \varepsilon$ and $|x|^p + |y|^p \leq K$ then $h(x | y) \geq \delta$. Indeed, suppose towards a contradiction that there exist $\varepsilon > 0$, $K > 1$ and sequences $\{x_n\}, \{y_n\} \subseteq \R^m$ such that $|x_n - y_n| \geq \varepsilon$, $|x_n|^p + |y_n|^p \leq K$, and $h(x_n | y_n) < 1/n$ for all $n \in \N$. Then $\{x_n \}, \{ y_n\}$ are bounded sequences, therefore there are subsequences $\{x_{k_n} \}, \{ y_{k_n}\}$ that converge to $x,y \in \R^m$, respectively. Note that $x \neq y$, since $|x - y| \geq \varepsilon$, and $h(x | y) = 0$ by continuity, which contradicts the strict convexity of $h$, thus establishing the claim. \par 
Let $\varepsilon > 0$, choose any $K > 1$, and let $\delta > 0$ be as in the previous claim. Then 
\begin{align*}
|\{ |u_n - u| \geq \varepsilon \}|  \leq \ &   |\{ |u_n - u| \geq \varepsilon \} \cap \{ |u_n|^p + |v_n|^p \leq K \}| \\
& + |\{ |u_n - u| \geq \varepsilon \} \cap \{ |u_n|^p + |v_n|^p > K \}| \\
 \leq \ &  |\{ h(u_n, u) \geq \delta  \}| + |\{ |u_n|^p + |v_n|^p > K \}| \\
 \leq \ & \frac{1}{\delta} \int_\Omega h(u_n | u) \, \rd x + \frac{1}{K} \int_\Omega |u_n|^p + |u|^p \, \rd x
\end{align*}
which implies that 
\[\lim_{n \to \infty} |\{ |u_n - u| \geq \varepsilon \}| \leq \frac{C}{K}\]
for some $C>0$. Since $\varepsilon$ and $K$ are arbitrary, letting $K \to \infty$ yields that $u_n$ converges to $u$ in measure, thereby completing the proof.

\section{Auxiliary classical results} \label{appendix_B}
\begin{definition}
Let $p \geq 1$ and $(X, \Sigma, \mu)$ be a measure space. A sequence $\{u_n\}$ in $L^p(X)$ is said to be \textit{uniformly $p$-integrable} if for every $\varepsilon>0$ there exists $\delta_\varepsilon>0$ such that 
\[\sup_n \int_E |u_n|^p \, \rd \mu < \varepsilon^p \]
for all $E\in \Sigma$ with $\mu(E) < \delta_\varepsilon$.
A sequence that is uniformly $1$-integrable is simply called \textit{uniformly integrable}.
\end{definition}

\begin{proposition}\label{prop:uniform_integrability_characterization}
Let $p \geq 1$ and $(X, \Sigma, \mu)$ be a finite measure space. If $\{u_n\}$ is a bounded sequence in $L^p(X)$ then $\{u_n\}$ is uniformly $p$-integrable if, and only if, for every $\varepsilon > 0$ there exists $M_\varepsilon > 0$ such that 
\[\sup_n \int_{ \{|u_n| \geq M_\varepsilon \}  } |u_n|^p \, \rd \mu < \varepsilon^p. \]
\end{proposition}
\noindent The proof of Proposition \ref{prop:uniform_integrability_characterization} for $p=1$ can be found in \cite[Chapter 19.5]{royden2010real}. For $p>1$ it extends easily by noting that $\{u_n\}$ is uniformly $p$-integrable if and only if $\{|u_n|^p\}$ is uniformly integrable.

\begin{theorem}[Vitali Convergence Theorem \cite{bartle1995elements}] \label{vitalithm}
Let $p \geq 1$ and $(X, \Sigma, \mu)$ be a finite measure space. Consider a sequence $\{u_n\}$ in $L^p(X)$ and a measurable function $u$. Then $\{u_n\}$ converges to $u$ in $L^p(X)$ if, and only if, $\{u_n\}$ is uniformly $p$-integrable and converges to $u$ in measure.
\end{theorem}

\begin{theorem}[Dunford-Pettis Theorem \cite{royden2010real}] 
 Let $(X, \Sigma, \mu)$ be a finite measure space. If $\{u_n \}$ is a bounded sequence in $L^1(X)$, then $\{u_n \}$ is relatively weakly compact if, and only if, $\{u_n \}$ is uniformly integrable.
\end{theorem}

\begin{proposition} \label{prop_appendix}
    Let $p \geq 1$ and $\Omega$ be a bounded open subset of $\R^d$, $d \in \N$. Consider a sequence $\{ u_n\}$ of functions $u_n : \Omega \to \R^n$, $n \in \N$, and assume that $\{u_n \}$ is bounded in $L^p(\Omega)$. Let $f : \R^n \to \R$ be a continuous function satisfying \[ \frac{|f(u)|}{1+|u|^p} \to 0 \qquad \mbox{as $|u| \to \infty$.} \]
    Then, the sequence $\{f(u_n) \}$ is relatively weakly compact in $L^1(\Omega)$.
\end{proposition}

\end{document}